\documentclass[letterpaper, 10 pt, conference]{ieeeconf}  

\IEEEoverridecommandlockouts                              
\title{\LARGE \bf
Decoupling Collision Avoidance in and for Optimal Control using Least-Squares Support Vector Machines
}

\author{Dries Dirckx and Wilm Decr\'{e} and Jan Swevers
\thanks{This work is funded by Flanders Make through the SBO project ARENA: Agile and
Reliable Navigation.}
\thanks{D. Dirckx, W. Decr\'{e} and J. Swevers are with Department of Mechanical Engineering, Faculty of Engineering Science,
        KU Leuven, 3001 Heverlee, Belgium
        {\tt\small dries.dirckx@kuleuven.be}}%
}

\usepackage{textcomp} 
\usepackage{pifont}   
\usepackage{booktabs}
\usepackage{amssymb, amsfonts} 
\usepackage{amsmath}
\usepackage{multicol}
\usepackage{graphicx}
\usepackage{textcomp}
\usepackage{xcolor}
\usepackage[utf8]{inputenc}
\usepackage{pgfplots}
\usepackage{pdfpages}
\usepackage{tabularx}
\usepackage{tikz}
\usepackage[font=normalsize]{caption}
\usepackage{tgadventor}
\usetikzlibrary{positioning}
\usepackage{svg}
\usepackage{arydshln}
\usepackage{tcolorbox}
\tcbuselibrary{theorems}
\usepackage{algorithm2e}
\usepackage{gensymb}
\usepackage{tikz-cd}

\definecolor{lightgraybox}{rgb}{0.85, 0.85, 0.85}

\newcolumntype{Y}{>{\centering\arraybackslash}X}

\usepackage[hyperref=auto, mincrossrefs=999, backend=biber, style=numeric-comp]{biblatex}
\usepackage{mathtools}
\usepackage{xurl}  

\usepackage{makecell}
\usepackage{multirow}
\usepackage{optidef}
\usepackage{tikz-cd}

\renewcommand{\vector}[1]{\boldsymbol{#1}}

\makeatletter
\newcommand*\mysize{%
  \@setfontsize\mysize{7}{14}%
}
\makeatother

\usepackage[english]{babel}
\usepackage{amsthm}
\newtheorem{theorem}{Theorem}
\theoremstyle{definition}
\newtheorem{problem}{Problem}

\newcolumntype{Y}{>{\centering\arraybackslash}X}
\newcolumntype{P}[1]{>{\centering\arraybackslash}p{#1}}

\newcommand\Set[2]{\{\,#1\mid#2\,\}}

\usepackage{array}
\usepackage{stfloats}
\usepackage{float}
\usepackage{verbatim}
\usepackage[T1]{fontenc}
\usepackage{subcaption}
\usetikzlibrary{shapes}
\usepackage[export]{adjustbox}
\usepackage{contour}
\usepackage{rotating}
\usepackage[normalem]{ulem}
\newcommand{\green}[1]{{\color{green}#1}}

\renewcommand{\vec}[1]{\boldsymbol{#1}}

\DeclareRobustCommand\sampleline[1]{%
  \tikz\draw[#1] (0,0) (0,\the\dimexpr\fontdimen22\textfont2)
  -- (1.5em,\the\dimexpr\fontdimen22\textfont2);%
}

\newcommand{\splitline}[4]{
\tikz{\draw[#1, #3, #4] (0, 0) (0,  \the\dimexpr\fontdimen22\textfont2) -- (0.75em,  \the\dimexpr\fontdimen22\textfont2) edge[#2] +(0:0.75em);}
}

\newtcbtheorem[number within=section]{mytheo}{}%
{colback=lightgray!15,colframe=gray!85, separator sign none}{th}

\definecolor{green}{RGB}{0, 170, 0}
\definecolor{red}{RGB}{200, 0, 0}
\definecolor{lightgray}{RGB}{200, 200, 200}
\definecolor{darkgray}{RGB}{150, 150, 150}
\definecolor{lightblue}{RGB}{41, 200, 200}
\definecolor{darkblue}{RGB}{41, 92, 207}
\definecolor{jump_red}{RGB}{176, 33, 33}
\definecolor{jump_yellow}{RGB}{250, 191, 5}
\definecolor{steelblue}{RGB}{70, 130, 180}
\definecolor{firebrick}{RGB}{178, 34, 34}

\contourlength{0.8pt}

\newcolumntype{L}[1]{>{\raggedright\let\newline\\\arraybackslash\hspace{0pt}}p{#1}}
\newcolumntype{C}[1]{>{\centering\let\newline\\\arraybackslash\hspace{0pt}}p{#1}}
\newcolumntype{R}[1]{>{\raggedleft\let\newline\\\arraybackslash\hspace{0pt}}p{#1}}

\begin{document}

\maketitle
\thispagestyle{empty}
\pagestyle{empty}

\begin{abstract}

This paper details an approach to linearise differentiable but non-convex collision avoidance constraints tailored to convex shapes. It revisits introducing differential collision avoidance constraints for convex objects into an optimal control problem (OCP) using the separating hyperplane theorem. By framing this theorem as a classification problem, the hyperplanes are eliminated as optimisation variables from the OCP. This effectively transforms non-convex constraints into linear constraints. A bi-level algorithm computes the hyperplanes between the iterations of an optimisation solver and subsequently embeds them as parameters into the OCP. Experiments demonstrate the approach's favourable scalability towards cluttered environments and its applicability to various motion planning approaches. It decreases trajectory computation times between 50\% and 90\% compared to a state-of-the-art approach that directly includes the hyperplanes as variables in the optimal control problem.

\end{abstract}

\section{Introduction: Autonomous Motion Planning}
\label{sec:introduction}
Deploying autonomous robots in practical and real-life settings, e.g., a warehouse, industrial production cell, homes, etc. is a complex problem with many interesting challenges that remain. These robots need to be aware of and interact with various actors, e.g., other robots, humans, their (immediate) surroundings, etc. Therefore, guaranteeing safe placing of autonomous robots in these settings continues to be one of the main pillars of new control algorithms, which possibly has severe consequences if safety considerations are not properly taken into account. These algorithms need to incorporate, among others, human-interpretable behaviour and operational (safety) guarantees, such as preventing collisions with humans, workspace objects, walls, tables, etc. and must adhere to the robot geometrical and kinodynamical restrictions such as the reachable workspace and actuator limits.\\

Current developments in automation and robotics increasingly move towards quick and (semi-) autonomous re-deployment of robot planners in new environments and contexts to allow for tailored small-batch production processes instead of mass production. This transition requires generic approaches for collision avoidance that are not only agnostic to a specific type of manipulator and environment but also maintain computational efficiency and differentiability. Trajectory planners simultaneously generate a geometric path and augment it with timing information and corresponding actuator input, e.g., velocity, acceleration or torque, to track that path, often called a motion profile. A subset of optimisation-based trajectory planning, called optimal control-based methods, solve an optimal control problem (OCP) to find such a trajectory~\cite{diehl2011numerical}. An OCP is formulated in continuous time and since finding an analytical solution is often not possible or tractable, it is often discretised to a nonlinear problem (NLP) and numerically solved~\cite{grandia2019, Foehn2021, Verschueren2014}. Popular mature and robust numerical solving techniques include sequential quadratic programming~\cite{acados}, interior-point methods~\cite{Waechter2005}, etc. These solvers typically rely on the first- and second-order derivatives of the constraints in the nonlinear programming (NLP) problem, thereby requiring differentiable collision detection.\\ 

The work on this paper revisits an existing differentiable, generally applicable collision detection constraint for optimal control-based trajectory planning. It introduces a reformulation of the existing non-convex constraint based on separating hyperplanes into a linear constraint at the cost of solving an extra linear system or quadratic program. Furthermore, it decouples the collision avoidance constraint from the OCP which enables to reduce the computational complexity of solving the OCP by removing unnecessary variables and constraints from the OCP.

\section{Related Work}
\label{sec:related_work}

\subsubsection*{Computational Geometrics}
Currently, many academic and industrial robot planners represent geometrical information of the environment either through meshes (typically in sampling-based planners) or as convex hull approximations, e.g., cuboids, capsules, etc. (mainly in optimisation-based planners). Many existing, computationally efficient toolboxes, often developed for fast collision detection, computer graphics or physics simulators, use geometric information to efficiently check if two convex objects are in contact~\cite{gjk1988, Montaut2022}. They are based on the famous Gilbert-Johnson-Keerthi algorithm and Expanding Polytope Algorithm for finding contact and separation vectors between convex shapes. The combination of both methods is popular due to their numerical robustness, high efficiency and ability to deal with generic shapes such as meshes. However, the functions provided by these toolboxes are not continuously differentiable and as such, not directly usable in optimisation-based techniques for trajectory optimisation that rely on derivative-based solvers.\\

\subsubsection*{Optimisation-based Collision Detection}
Recently, several efforts have been made to create differentiable formulations for collision detection. Randomised smoothing, from the field of reinforcement learning, presents an interesting approach to construct derivatives based on multiple evaluations of these geometric detection toolboxes. Montaut et al.~\cite{Montaut2023} use this concept to efficiently find local contact derivative information of geometric meshes used in differentiable physics simulations, using a performant implementation of the GJK algorithm. Nonetheless, this method generates stochastic gradients and is not (yet) applicable to optimisation-based trajectory planning using gradient-descent solvers. \\


Alternatively, Tracy et al.~\cite{Tracy2023} formulate a convex optimisation problem with conic constraints to detect collision between convex primitives, using the implicit function theorem to calculate derivatives of that convex problem to the convex primitives involved in collision checking. However, solving a conic program requires multiple solutions to a linear system (the KKT-system) which are quite costly to obtain. As collision checks between every pair of robot-object primitives in the environment must be performed several times within a planner, this quickly leads to a significant amount of collision checks.\\


\subsubsection*{Differentiable Collision Detection}
The idea of considering collision detection as a classification problem has already been used before in literature. Zhi et al.~\cite{Zhi2022} map the robot's configuration space to (multi-)class labels using differentiable kernels based on a robot's forward kinematics. Firstly, their classifier update procedure heavily relies on the assumption that the obstacles do not move too fast and thus cannot be used for a new obstacle or environment. Secondly, the classifier is trained for a specific robot geometry and requires training a new classifier for a different robot (geometry) or environment. Mercy et al.~\cite{Mercy2017} and Vandewal et al.~\cite{Vandewal2020} use the separating hyperplane theorem~\cite{Boyd2004} in spline-based motion planning due to its polynomial structure. It is used to guarantee continuous collision avoidance by approximating state trajectories with B-splines or Bernstein polynomials. They incorporate the hyperplanes directly as optimisation variables within the optimal control problem (OCP), resulting in bilinear, non-convex constraints that makes the OCP harder to solve.


\section{Contributions}
In this work, the geometric information of the robot and obstacle in Cartesian space is considered. We propose to formulate collision detection with the separating hyperplane theorem but not embed the hyperplane parameters as optimisation variables. By considering the robot and obstacle as separate labeled classes, the hyperplanes are found separately as the solution of a least-squares classification problem, creating a linearised version of a bilinear constraint~\cite{Mercy2017, Vandewal2020}. This work contributes by:
\begin{enumerate}
    \item formulating collision avoidance for general convex shapes as a least-squares classification problem,
    \item casting the non-convex separating hyperplane constraint into a deterministically differentiable, generic, linear constraint,
    \item creating a robot-agnostic collision avoidance constraint based on classification techniques.
\end{enumerate}

The paper is structured as follows. Firstly, Section~\ref{sec:lssvm} covers the formulation of the supervised classification problem as a quadratic program (QP-SVM) or a linear system (LS-SVM). Secondly, Section~\ref{sec:bilevel} covers the integration of the hyperplanes as parameters in an optimisation problem for trajectory planning, leading to a two-step algorithm. Lastly, the presented reformulation's performance is compared to the approach in~\cite{Mercy2017} on two experimental cases: motion planning for a holonomic vehicle and a 6-DOF robotic manipulator.

\section{Collision Detection as Least-Squares Classification}
\label{sec:lssvm}
Consider the robot and obstacles fully described by (a combination of) spheres, capsules or convex polytopic primitive(s). The robot is modelled with $N_{r}$ primitives $V_{r}$, while all $m$ obstacles are jointly described with $N_{o}$ primitives $V_{o}$. The set $R = \Set{v_{i}}{v_{i} \in V_{r}}$ and $O = \Set{v_{k}}{v_{k} \in V_{o}}$ contain all vertices of the robot and obstacles, respectively. The vertices $v \in \mathbb{R}^{2},\;\mathbb{R}^{3}$ for a single primitive can be seen as data points belonging to a specific class, with one primitive assigned to a single class. Two such convex sets A and B are separable by a hyperplane with normal $\boldsymbol{w}$ and offset $w_{b}$ under mild conditions~\cite{Boyd2004}: 

\begin{theorem}[Separating hyperplane theorem~\cite{Boyd2004}]
Suppose $A$ and $B$ are nonempty disjoint convex sets, i.e., $A \cap B$ = $\emptyset$. Then there exist $\boldsymbol{w} \neq 0$ and $w_{b}$ such that $\boldsymbol{w}^{\top}\boldsymbol{y} + w_{b} \geq 0$ for all $\boldsymbol{y} \in A$ and $\boldsymbol{w}^{\top}\boldsymbol{y} + w_{b} \leq 0$ for all $\boldsymbol{y} \in B$. The hyperplane $\boldsymbol{w}^{\top}\boldsymbol{y} + w_{b} = 0$ is called a \textit{separating hyperplane} for the sets A and B.
\label{theorem:sep_hyperplanes}
\end{theorem}

\begin{figure}[h]
    \begin{minipage}[c]{0.5\linewidth}
        \centering
        \includegraphics[width=.95\textwidth]{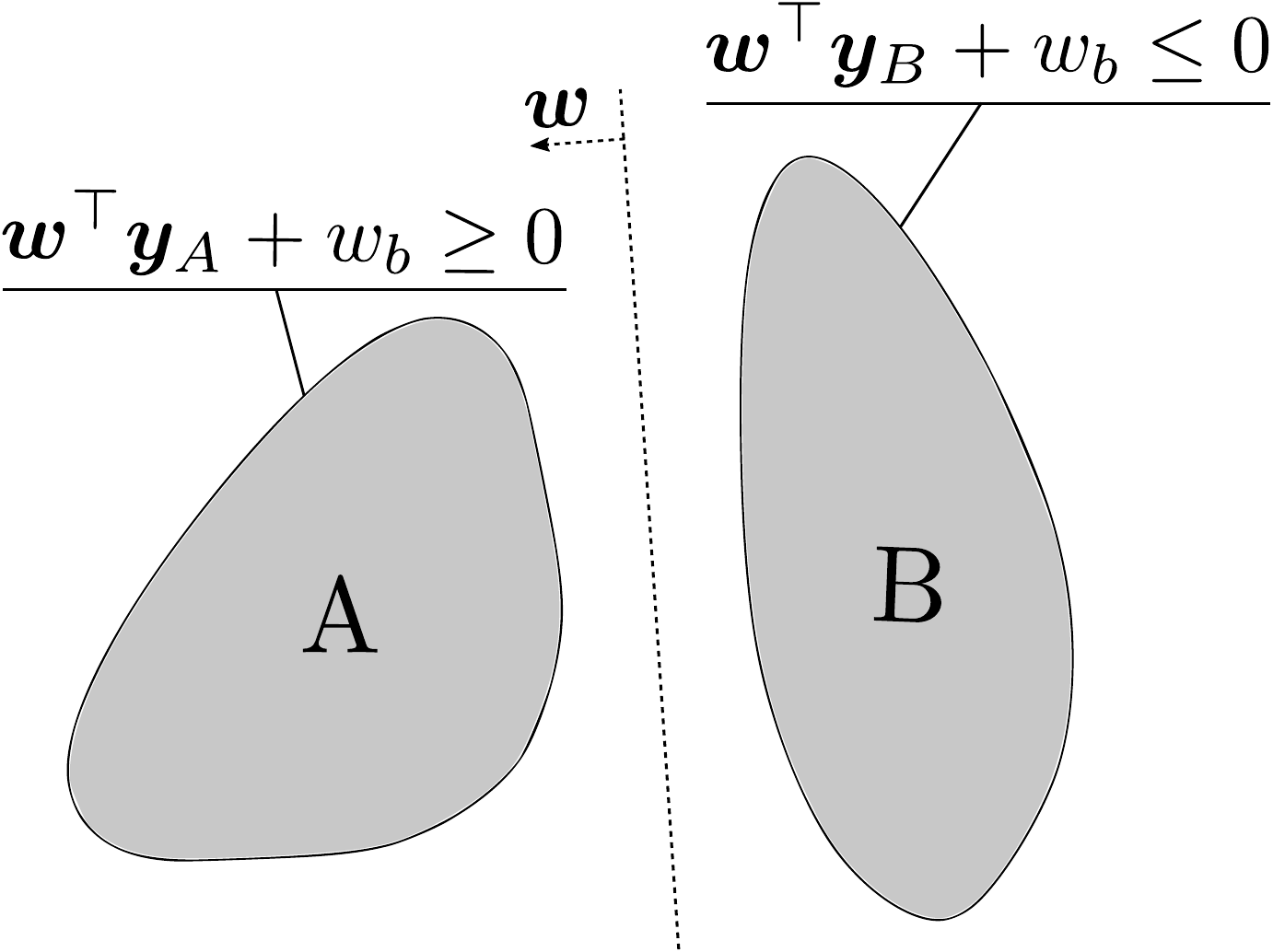}
    \end{minipage}
    \hspace{-3ex}
    \begin{minipage}[c]{0.5\linewidth}
        \centering
        \includegraphics[width=0.95\textwidth]{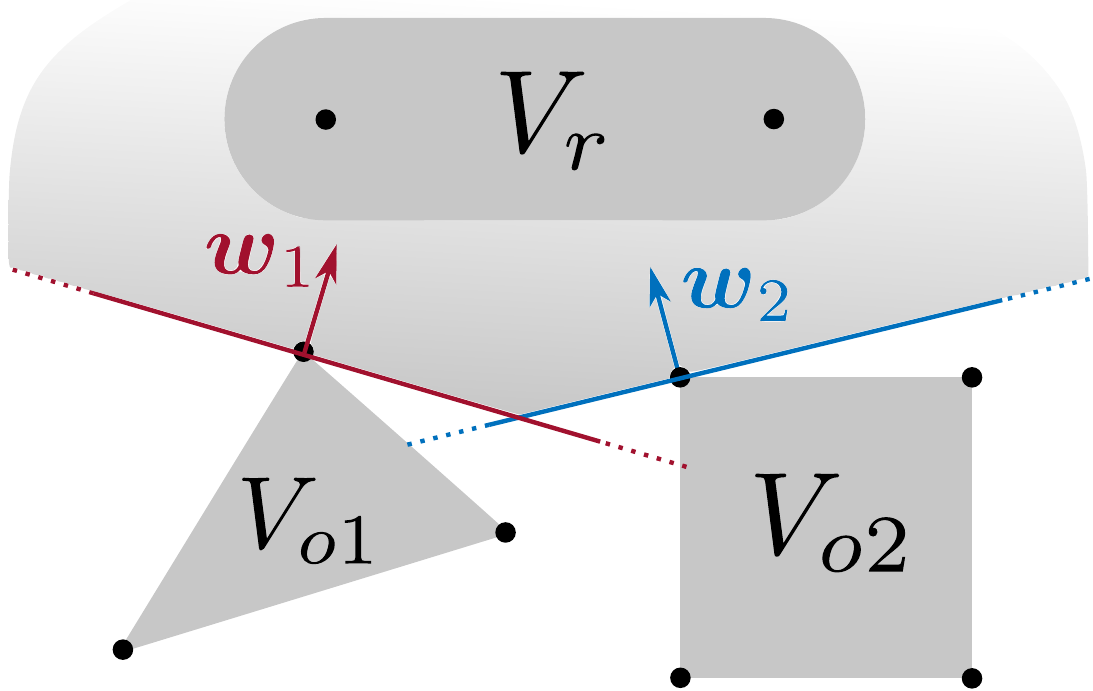}
    \end{minipage}
    \caption{(left) Separating hyperplane for two disjoint, nonempty convex sets A and B~\cite{Boyd2004}. (right) Illustration of the set of vertices for a capsular robot ($V_{r}$) and two obstacles, one triangular ($V_{o1}$) and one rectangular ($V_{o2}$).}
    \label{fig:hyper_theorem}
\end{figure}
Figure~\ref{fig:hyper_theorem} illustrates the theorem for general sets and for three geometric primitives: a two-dimensional capsule, a rectangle, and a triangle. The robot $V_{r}$ must always remain in the (half-open) polygonal approximation of the free space, constructed by the hyperplanes $\boldsymbol{w}_{1}$ and $\boldsymbol{w}_{2}$. The content in this paper focuses on robots and obstacles described geometrically in Cartesian space. In this case, a hyperplane that satisfies Theorem~\ref{theorem:sep_hyperplanes} corresponds to a linear classification boundary in Cartesian space. The remainder of this work focuses on linearly separable classes.





\subsection{(Least-Squares) Support Vector Machines}
Sets A and B consist of discrete labeled data points and finding the separating hyperplane corresponds to finding a solution to a supervised classification problem. This problem can be solved with the technique of Support Vector Machines~\cite{Vapnik1995}, formulating the supervised classification problem as a quadratic program (QP). 

\begin{problem}[Support Vector Machines (QP-SVM)~\cite{Vapnik1995}]{}
Suppose a set of $N_{d}$ training data points $\text{\{} \Gamma_k, y_k\text{\}}^{N_{d}}_{k=1}$, where $y_{k}$ $\in \mathbb{R}^{n_{y}}$ are the data points and $\Gamma_k$ the corresponding class labels $\text{\{}-1, 1\text{\}}$. The maximal-margin linear classification boundary, characterised by vector $\boldsymbol{w}$ and offset $w_b$, is represented as the solution to the following QP:
\begin{equation}
    \begin{aligned}
      \underset{\boldsymbol{w}, w_{b}, e_{k}}{\textbf{minimise}} \quad & \frac{1}{2} \boldsymbol{w}^{\top}\boldsymbol{w} + \tau\sum_{k=1}^{N_d} e_{k} \\
      \textbf{subject to} \quad & \Gamma_{k}[\boldsymbol{w}^{\top}\Xi({\boldsymbol{y}_{k}}) + w_b] \geq 1 - e_{k},\\
      & e_{k} \geq 0.
    \end{aligned}
    \label{eq:qpsvm}
\end{equation}
\end{problem}
The variable $e_{k}$ represents a possible misclassification, visualised in Fig.~\ref{fig:linear_svm}. $\Xi(\cdot)$ maps the data $\boldsymbol{y}_{k}$ to a higher dimensional space (the `kernel trick') to deal with nonlinearly separable sets and $\tau > 0$ penalises the misclassification $e_k$. This penalisation allows us to have a larger margin $\frac{1}{\boldsymbol{w}^{\top}\boldsymbol{w}}$ around the classification boundary, to achieve a clearer separation of the classes, at the cost of misclassifying a few data points. For linearly separable classes, as considered in this work, the error terms $e_{k}$ are not necessary. They are ommited from the QP which is the QP that is used in Section~\ref{sec:experiments}. The LS-SVM approach from Suykens et al.~\cite{Suykens1999} squares the penalty term in the objective function of Problem~\ref{eq:qpsvm} and reformulates the inequality constraints as equalities, resulting in a convex QP with equality constraints. As this reformulation satisfies Slater's condition, there is no duality gap and the solution to their problem can be found as the solution to the dual problem. The dual problem is formulated as the following dense system of linear equations:

\begin{figure}
    \centering
    \includegraphics[width=0.5\linewidth]{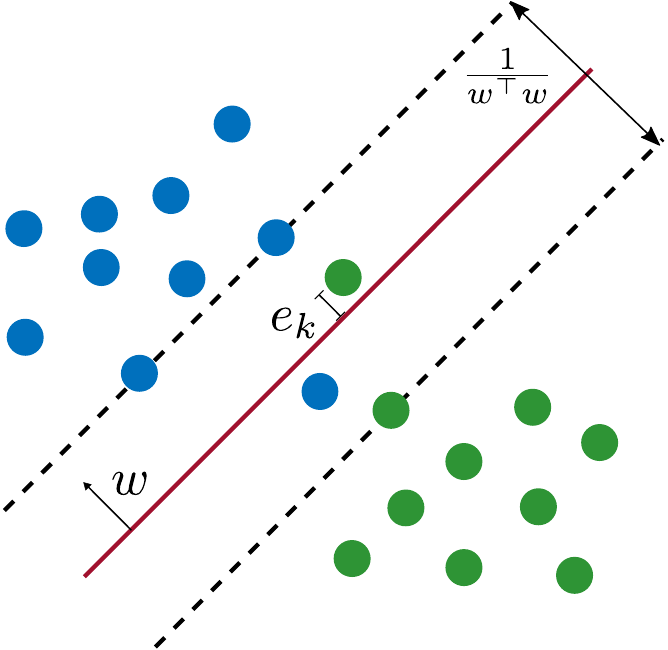}
    \caption{Illustration of two linearly separable classes, the hyperplane normal $w$, the classification margin $\frac{1}{w^{\top}w}$ and a misclassification characterised by the error $e_{k}$.}
    \label{fig:linear_svm}
\end{figure}



\begingroup
\setlength{\arraycolsep}{1pt} 
\begin{equation}
    \begin{aligned}
        \begin{bmatrix}
            0 & -\Gamma_{1} & \ldots & -\Gamma_{N_{d}} \\[0.5em]
            \Gamma_{1} & \boldsymbol{Z}_{1*}\boldsymbol{Z}_{1*}^{\top} + \frac{1}{\tau} & \ldots & \boldsymbol{Z}_{1*}\boldsymbol{Z}_{N_{d}*}^{\top} \\
            \vdots & \vdots & \ddots & \vdots \\
            \Gamma_{N_{d}} & \boldsymbol{Z}_{N_{d}*}\boldsymbol{Z}_{1*}^{\top} & \ldots & \boldsymbol{Z}_{N_{d}*}\boldsymbol{Z}_{N_{d}*}^{\top} + \frac{1}{\tau}
        \end{bmatrix}
        \begin{bmatrix}
            w_{b} \\
            \boldsymbol{\alpha}
        \end{bmatrix} = 
        \begin{bmatrix}
            0 \\
            \vec{1}
        \end{bmatrix}
    \end{aligned}
    \label{eq:lssvm}
\end{equation}
\endgroup

where $\boldsymbol{Z} \in \mathbb{R}^{N_d \times n_{y}}$ is the matrix of labeled data points $[\Xi(\boldsymbol{y}_{1})^{\top}\Gamma_1,\, \ldots,\, \Xi(\boldsymbol{y}_{N_d})^{\top}\Gamma_{N_{d}}]$, $\vec{1} \in \mathbb{R}^{N_d}$ and $\boldsymbol{\alpha} = [\alpha_{1}, \ldots, \alpha_{N_d}]$ represent the Lagrangian multipliers of the LS-SVM problem. By eliminating $\boldsymbol{\alpha}$, the Cholesky refactorisation can be applied to the Schur complement of the above matrix to solve Eq.~\ref{eq:lssvm} and allows us to find $\boldsymbol{w} = \sum_{k=1}^{N_{d}} \Gamma_k \alpha_k \boldsymbol{y}_k$ and $w_{b}$. Solving Eq.~\ref{eq:qpsvm} reduces the conservativeness introduced by the dense solution of Eq.~\ref{eq:lssvm} at the cost of a higher computational effort.\\


Since this work considers convex objects that should be separable with a hyperplane in Cartesian space, the kernel function $\Xi$ is chosen as the linear kernel $\Xi(\boldsymbol{y}_{i}) = \boldsymbol{y}_{i}$. Penetration of robot and object is accommodated in Eq.~\ref{eq:lssvm} with the misclassification term $e_{k}$. Allowing misclassification leads to a hyperplanes that can still be defined if the robot and obstacles are in contact. This allows us to start from an infeasible initial guess of the robot's trajectory.


\section{Decoupled Collision Avoidance in Optimal Control}
\label{sec:bilevel}
Using Theorem~\ref{theorem:sep_hyperplanes} for collision avoidance, a discrete model for the system $\vec{x}_{k + 1}= F(\vec{x}_{k}, \vec{u}_{k}, T)$ in terms of the system's states $\vec{x}_{k}$ and controls $\vec{u}_{k}$, and a trajectory cost $\ell$ leads to the following general optimal control problem for point-to-point motion planning: 
\begin{equation}
    \begin{aligned}{}
        \centering
        \underset{\substack{\vector{x}_{k}, \vector{u}_{k}, T, \\ 
        \vec{w}_{s}, \vec{w}_{b, s}}}{\textbf{minimise}} \quad & \sum_{k=0}^{N - 1}\ell(\vec{x}_{k}, \vec{u}_{k}, T) + \ell_{N}(\vec{x}_{N})\\
        \textbf{subject to} \quad & \vec{x}_{k + 1}= F(\vec{x}_{k}, \vec{u}_{k}, T), \\
        & \boldsymbol{w}_{s}^{\top}\boldsymbol{p}_{r} + w_{b, s} \geq \mathcal{D}, \\
        & \boldsymbol{w}_{s}^{\top}\boldsymbol{p}_{o} + w_{b, s} \leq 0, \\
        & \| \boldsymbol{w}_{s} \|^2 = 1,  \\
        & \vec{x}_{0}= \vec{x}_{\text{start}}, \; \vec{x}_{N} = \vec{x}_{\text{goal}}, \\
        & \vec{x}_{k}, \vec{u}_{k} \in \mathbb{X}, \mathbb{U}, \\
        & s = 1, \ldots, N_{r}N_{o} \\
        & k = 0, \ldots, N - 1
        \label{eq:prelim_nlp_motplan}
    \end{aligned}
\end{equation}
where T is the total trajectory time, $\mathbb{X}, \mathbb{U}$ the set of feasible states and controls, $\boldsymbol{p}_{r}$ and $\boldsymbol{p}_{o}$ represent each vertex's Cartesian position for the robot and obstacles in $V_{r}$ and $V_{o}$ and $\mathcal{D}$ represents the relevant dimension of the robot and/or obstacle, i.e., the radius for a sphere or capsule primitive. Collision avoidance between each robot-obstacle pair is represented by the hyperplane normal $\boldsymbol{w}_{s}$ and offset $w_{b, s}$.\\



This work proposes to consider the hyperplanes optimisation variables $\boldsymbol{w}_{s}, w_{b, s}$ as optimisation problem parameters and update their values between iterations using Eq.~\ref{eq:qpsvm} and~\ref{eq:lssvm}. Consequently, this eliminates all optimisation variables from certain constraints as only the robot position remains as variable, marked as green below. These constraints are no longer enforceable in an optimisation problem and are removed, effectively reducing the amount of constraints:
\vspace{-0.25ex}
\begin{equation}
\left.
\begin{aligned}
            \boldsymbol{w}_{s}^{\top}\green{\boldsymbol{p}_{r}} + w_{b, s} &\geq \mathcal{D}\\
            \boldsymbol{w}_{s}^{\top}\boldsymbol{p}_{o} + w_{b, s} &\leq 0\\
            \| \boldsymbol{w}_{s} \|^2 &= 1
\end{aligned}
\quad
\right\}
\quad
\longrightarrow
\quad
\boldsymbol{w}_{s}^{\top}\green{\boldsymbol{p}_{r}} + w_b \geq \mathcal{D},
\label{eq:hyperplane_decoupled}
\end{equation}
which can also be found as the zero-order Taylor approximation of the multivariate function $g = \boldsymbol{w}_{s}^{\top}\boldsymbol{p}_{r} + w_{b, s}$. To enforce collision avoidance without an explicit constraint on the obstacle side, the hyperplane is shifted to the closest vertex $\boldsymbol{v}_{c}$ along its normal $\boldsymbol{w}$ using a new offset $w'_{b} = -\boldsymbol{w}^{T}\boldsymbol{v}_{c}$.

Removing the hyperplanes from the optimisation variables not only reduces the amount of optimisation variables but also eliminates a total amount of $m(N_{o} + 1)N$ constraints, with N the amount of control intervals of the OCP. This greatly reduces the cost of solving Eq.~\ref{eq:prelim_nlp_motplan} but requires solving Eq.~\ref{eq:lssvm} multiple times.

\subsection{Broad-phase and Trust-region Filter}
The hyperplanes are considered as optimisation problem parameters and are thus no longer computed iteratively as part of the current iterate $\boldsymbol{x}_{v}$. In the decoupled approach, their values $\boldsymbol{w}_{s}$ and $w_{b, s}$ would have to be updated each iteration to compensate for this fact. This means that the constraints of the NLP change during the solving process, which interferes with globalisation strategies of NLP solvers~\cite{nocedal, Waechter2005, snopt}. To limit the frequency at which these problem parameters change between iterations, and thus the negative effect on the globalisation, two filters are implemented: a broad-phase and a trust-region filter. \\

Firstly, a \textit{broad-phase filter} eliminates unnecessary updates of the hyperplanes for obstacles at a fixed distance $d_{bp}$ away from the robot. The effect is two-fold: (1) decreasing the amount of times Eq.~\ref{eq:qpsvm} or~\ref{eq:lssvm} is solved and (2) excluding unnecessary changes in the constraints between iterations. Secondly, to prevent numerical difficulties close to convergence due to small changes in hyperplanes over subsequent iterations, a simple \textit{trust-region strategy} is adopted. As a solver approaches a locally optimal trajectory, small iterative steps from $\boldsymbol{x}_{v - 1}$ to $\boldsymbol{x}_{v}$ are taken by the solver to make a stable approach towards that optimum. However, changing the hyperplane parameters in this phase can deteriorate convergence. At iteration $v$, the constraint for a point $\boldsymbol{p}_{i}$ on trajectory $\boldsymbol{x}_{v}$ equals $\boldsymbol{w}_{v}^{\top}\boldsymbol{p}_{i} + w_{b, v}$ (a). Subsequently, the solver computes a new trajectory with a lower constraint value equal to $\boldsymbol{w}_{v}^{\top}\boldsymbol{p}_{i, v + 1} + w_{b, v}$ (b). For this new trajectory, new hyperplanes $\boldsymbol{w}_{v + 1}$ are computed that differ from $\boldsymbol{w}_{v}$ (c). Consequently, it might happen that the new constraint $\boldsymbol{w}_{v + 1}^{\top}\boldsymbol{p}_{i, v + 1} + w_{b, v + 1}$ has a higher value than what the solver expected it to be, namely $\boldsymbol{w}_{v}^{\top}\boldsymbol{p}_{i, v + 1} + w_{b, v}$ (d). This would require additional iterations of the numerical solver to again reduce this constraint value and should be avoided. Therefore, if a hyperplane does not significantly change between two subsequent iterations $v$ and $v + 1$, the hyperplane of iteration $v$ is kept. The trade-off that is introduced with these filters is covered in Section~\ref{sec:experiments}.

\subsection{Two-step Optimisation}
The combination of (1) detaching the hyperplanes from a numerical solver and (2) employing filters to reduce negative effects of this detachment on convergence towards a local optimum, leads to the overall method detailed in Algorithms~\ref{alg:full_svm} to~\ref{alg:compute_svm}. The integration of the $\texttt{update}$ function into Ipopt is achieved through the Callback functionality of CasADi~\cite{Andersson2018}.

\RestyleAlgo{ruled}
\SetKwComment{Comment}{/* }{ */}
\begin{algorithm}[h]
    \caption{OCP with decoupled hyperplanes}\label{alg:full_svm}
    \KwData{Initial guess for state $\boldsymbol{x}_{\text{init}}$ and control $\boldsymbol{u}_{\text{init}}$}
    $\boldsymbol{w}_{0, *}, w_{b_{0, *}} \gets \texttt{update}(\boldsymbol{x}_{\text{init}}, R, O)$ \\
    \While{SOLVER not converged}{
    \tcc{One SOLVER iteration}
    $\boldsymbol{x}_{v}, \boldsymbol{u}_{v} \gets \texttt{OCP\_iter}
    (\boldsymbol{w}_{v-1, *}, w_{b_{v-1, *}})$ \\
    \tcc{One hyperplane update}
    $\boldsymbol{w}_{v, *}, w_{b_{v, *}} \gets \texttt{update}(\boldsymbol{x}_{v}, R, O, \boldsymbol{w}_{v-1}, w_{b_{v-1}})$
    }   
\end{algorithm}
\vspace{-4ex}
\begin{algorithm}[h]
    \caption{Inter-iteration hyperplane computation - $\texttt{update}$}\label{alg:update_svm}
    \KwData{Current solution $\boldsymbol{x}_{v}$, set $R$ and $O$, $\boldsymbol{w}_{v, *}, \, w_{b, v, *}$}
    $\texttt{x\_in\_coll} \gets \texttt{check\_trajectory}(\boldsymbol{x}_{v}, R, 0)$\\
    \eIf{$\texttt{x\_in\_coll}$}
    {
        \For{$i \in [0, N]$ \& $j \in [0,  N_{r}]$ \& $l \in [0, N_{o}]$}{
            \tcc{$* ={i, j, l}$}
            $\boldsymbol{w}_{v + 1, *} \gets \texttt{comp\_hyp}(V_{j} \, ; \, V_{l} \, ; \, \boldsymbol{w}_{v, *} \, ; \texttt{True})$ \\
            $w_{b, v + 1, *} \gets \texttt{vertex}(\boldsymbol{w}_{v + 1, *} \, ; \, V_{j} \, ; \, V_{l})$
        } 
    }{
        \For{$i \in [0, N]$ \& $j \in [0,  N_{r}]$ \& $l \in [0, N_{o}]$}{
            \tcc{$* ={i, j, l}$}
            $\boldsymbol{w}_{v + 1, *} \gets \texttt{comp\_hyp}(V_{j} ;  V_{l} ; \boldsymbol{w}_{v, *} ; \texttt{False})$\\
            $w_{b, v + 1, *} \gets \texttt{vertex}(\boldsymbol{w}_{v + 1, *}, V_{j}, V_{l})$
        } 
    }
\end{algorithm}
\vspace{-1ex}
\begin{algorithm}[h]
\caption{Computing the hyperplane normal $\boldsymbol{w}$ - $\texttt{compute\_hyp}$}\label{alg:compute_svm}
\KwData{Robot primitive $V_{j}$, obstacle primitive $V_{l}$, normal $\boldsymbol{w}_{v, i, j, l}$, $\texttt{LS}$}
\eIf{$\texttt{broad\_phase}(V_{j} \, ; \, V_{l}) \leq d_{bp}$}
{
    \tcc{Eq.~\ref{eq:qpsvm} or~\ref{eq:lssvm}}
    $\boldsymbol{w}'_{v + 1, i, j, l} \gets \texttt{solve\_svm}(V_{j} \, ; \, V_{l}  \, ; \, \texttt{LS})$
    $\theta_{tr} = \texttt{trust\_region}(\boldsymbol{w}_{v, i, j, l} \, ; \, \boldsymbol{w}'_{v + 1, i, j, l})$\;
    \If{$\theta_{tr} \leq \bar{\theta}$}{
        $\boldsymbol{w}_{v + 1, i, j, l} \gets \boldsymbol{w}'_{v + 1, i, j, l}$
    }
}
{$\boldsymbol{w}_{v + 1, i, j, l} \gets \boldsymbol{w}_{v, i, j, l}$}
\end{algorithm}

Algorithm~\ref{alg:full_svm} explains how the update of the hyperplanes is intertwined in the algorithm of a numerical solver. After an iteration, the update is performed and the hyperplane parameters $\boldsymbol{w}_{v}, w_{b, v}$ are used in the next iteration. Algorithm~\ref{alg:update_svm} details how for every point along the trajectory $\boldsymbol{x}_{v}$, a hyperplane and its corrected offset are computed for every pair of robot and obstacle primitives in the environment. A distinction is made between whether or not the trajectory $\boldsymbol{x}_{v}$ is collision-free, in which case Eq.~\ref{eq:qpsvm} without misclassification is solved to reduce conservativeness. If there still exists a point where there is a collision between any robot-obstacle, Eq.~\ref{eq:lssvm} is used to update the hyperplanes. Algorithm~\ref{alg:compute_svm} performs the filter checks per robot-obstacle primitive pair.
\vspace{-1.5ex}
\section{Experimental Validation}
\label{sec:experiments}


The presented method is studied in terms of various performance aspects on two cases: two-dimensional motion planning of a point mass and three-dimensional planning for a six-DOF robotic manipulator in a cluttered environment. For each case, a comparison is made between the decoupled problem formulation ('SVM') and a formulation where the hyperplanes are embedded as optimisation variables ('Coupled'). For each of the cases, parts of the initial trajectories will collide with one or more obstacles to investigate and show that the method can aptly handle partially infeasible initial guesses as fully feasible initialisations are often not available.\\

Table~\ref{tab:ocp_svm_details} lists all details necessary to fully describe the OCPs, referring back to the general form given in Eq.~\ref{eq:prelim_nlp_motplan}. Each OCP is discretised using multiple shooting with Runge-Kutta 4 as integrator into an NLP using CasADi. They are solved with Ipopt, using the linear solver \texttt{ma27}~\cite{HSL1983}. The linear solver ´lapacklu´ is used to solve the LS-SVM system of equations while ProxQP is used to solve Eq.~\ref{eq:qpsvm}. All experiments are conducted on a portable computer with an Intel® Core™ i7-10810U processor with twelve cores at 1.10GHz and 31.1 GB RAM.\\
\begin{table}[]
    \centering
    \addtolength{\tabcolsep}{8pt}
    \begin{tabular}{llll}
    \toprule
    \toprule
          & & \textbf{Holonomic} & \textbf{six-DOF robot} \\
    \midrule[0.05pt]
         \textbf{Objective} & $\ell$ & $\sum_{i}\|v_{i}\|^{2}$ & T + $\tau_{1}\|\ddot{\boldsymbol{q}}\|^{2}$\\
    \midrule[0.05pt]
         \textbf{States} & $\boldsymbol{x}$ & $\boldsymbol{p} = [p_{x}, p_y]$ & $[\boldsymbol{q}, \dot{\boldsymbol{q}}]$\\
    \midrule[0.05pt]
         \textbf{Controls} & $\boldsymbol{u}$ & $\boldsymbol{v} = [v_{x}, v_{y}]$ & $\ddot{\boldsymbol{q}}$  \\
    \midrule[0.05pt]
         \textbf{Dynamics} & f & $\frac{d}{dt}{\boldsymbol{p}} = \boldsymbol{v}$ & $\frac{d}{dt}{[\boldsymbol{q}, \dot{\boldsymbol{q}}]} = [\dot{\boldsymbol{q}}, \ddot{\boldsymbol{q}}]$\\
    \midrule[0.05pt]
         \textbf{Collision} & $g_{j}$ &
         \begin{tabular}{@{}l@{}}Eq.~\ref{eq:hyperplane_decoupled} \\ $\boldsymbol{p}_{r} \rightarrow \boldsymbol{p}$ \end{tabular} & 
         \begin{tabular}{@{}l@{}}Eq.~\ref{eq:hyperplane_decoupled} \\ $\boldsymbol{p}_{r} \rightarrow \kappa(\boldsymbol{q})$\end{tabular}\\
    \midrule[0.05pt]
        \textbf{\makecell[l]{Terminal \\ constraint}} & $g_{N}$ & $\boldsymbol{p}_{N} = \boldsymbol{p}_{goal}$ & $\kappa(\boldsymbol{q}_{N}) = \boldsymbol{p}_{f}$\\
    \midrule[0.05pt]
         \textbf{Parameters} & $d_{bp, 1}$ & [0.0, 0.15, 1.0] & 0.1\\
         & $d_{bp, 2}$ & [0.0, 0.15, 1.0] & 0.1\\
         & $\theta_{tr}$ & $[0.1\degree, 5\degree, 45\degree]$ & 5\degree\\
         & N & 30 & 21\\
         & $m$ & $[1, \cdots, 10]$ & 5\\
    \bottomrule
    \bottomrule
    \end{tabular}
    \caption{Overview of the formulations, chosen transcription methods of the OCP and parameters for the decoupled approach for each of the three cases.}
    \label{tab:ocp_svm_details}
\end{table}

\begin{figure}
    \centering
    \includegraphics[width=0.5\linewidth]{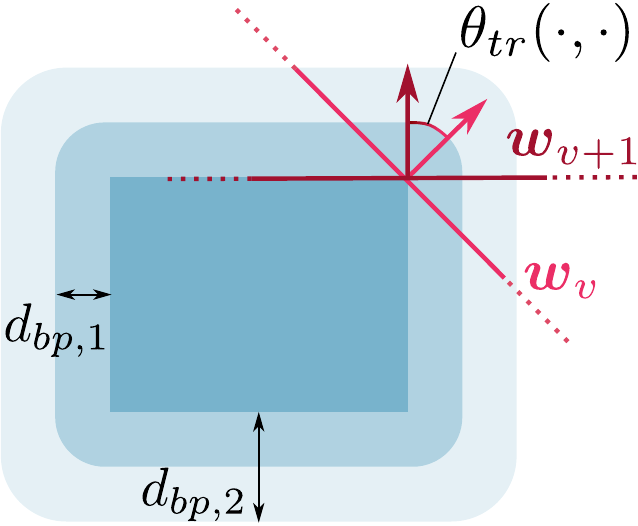}
    \caption{Visualisation of an example of the broadphase filter thresholds and the trust-region filter
described in Algorithm 3 and employed in the experiments.}
    \label{fig:filters}
\end{figure}

The broad-phase filter is implemented using a geometric distance function between a circle and a rectangle for the two-dimensional case and using the the geometric collision checker Coal~\cite{hppfcl} in the three-dimensional case. For the trust-region filter, the inter-normal angle $\theta_{tr} = \cos^{-1}(\frac{\boldsymbol{w}_{v}\cdot\boldsymbol{w}'_{v + 1}}{\|\boldsymbol{w}_{v}\| \|\boldsymbol{w}_{v + 1}\|})$ between two hyperplane normals is chosen in both cases. Both filter principles are visualised in Fig.~\ref{fig:filters}.


\subsection{Two-dimensional Holonomic Motion Planning}
\label{sec:holonomic}
\begin{figure}
    \begin{minipage}[c]{0.475\linewidth}
        \centering
        \includegraphics[width=0.8\textwidth]{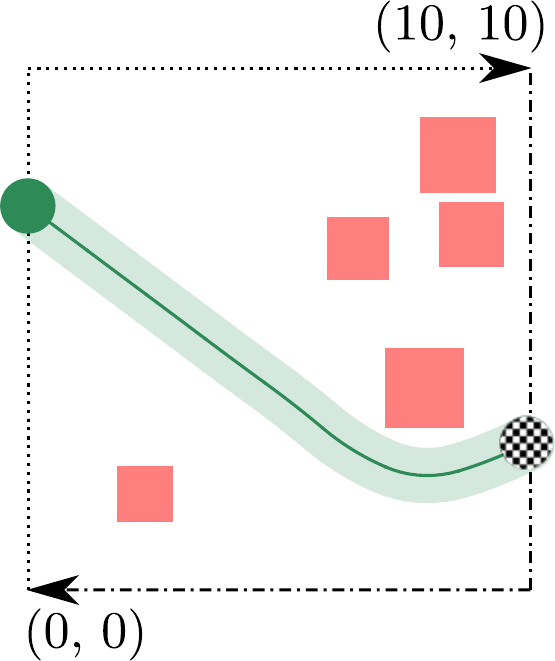}
    \end{minipage}
    \begin{minipage}[c]{0.475\linewidth}
        \centering
        \includegraphics[width=0.8\textwidth]{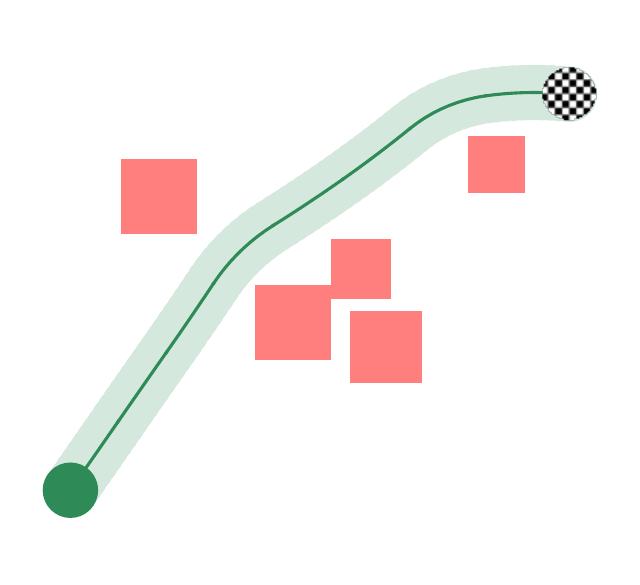}
    \end{minipage}
    \caption{Two randomly generated environments with five obstacles, where the AMR has to navigate through both in the holonomic case and spline-based case. The start and goal positions are indicated by the green and checkered circles, respectively. Ten different start and goal combinations are constructed per environment by varying them along the dotted and dash-dotted line.}
    \label{fig:obstacle_spline}
\end{figure}

A two-dimensional holonomic robot with circular footprint has to avoid rectangular obstacles, as illustrated in Fig.~\ref{fig:obstacle_spline}. The experiments are conducted for an increasing amount of obstacles in an environment, from one to ten obstacles. Per number of obstacles, a set of twenty environments is generated where the size of the obstacles is pseudo-randomly (uniform distribution) chosen between 0.3 and 0.8. For each environment, trajectories between ten start and end positions in the interval $[0, 10]$ must be computed. This leads to a total of 200 trajectories per number of obstacles. To ensure that the experiments do not just cover cases close to each other, the ten  combinations are chosen such that they cover all directions in which one could traverse an environment. The start position is varied from bottom left to top right while the goal position varies in the opposite direction. Each NLP is initialised with a straight line between start and end configuration, possibly colliding with obstacles in the environment. An initial guess for the hyperplane parameters is found by solving Eq.~\ref{eq:lssvm} for this straight-line trajectory.\\

\subsubsection*{Filter Effects}
Firstly, the effect of the filter parameters on the performance of the decoupled approach is studied. Table~\ref{tab:effects_filter} lists the results for the following validation metrics, calculated over 200 trajectories for an environment with four obstacles:
\begin{itemize}
    \item \# QP, \# LS: the average number of times Eq.~\ref{eq:qpsvm} and~\ref{eq:lssvm} are solved per trajectory
    \item \# I: the median amount of iterations required by Ipopt to solve the NLP
    \item Wall/I : the median wall time required per iteration of Ipopt, including the time needed to solve Eq.~\ref{eq:qpsvm} or~\ref{eq:lssvm}
    \item Opt. : the added optimality (shortest path) of the decoupled approach relative to the optimality of the coupled approach, calculated as the percentage $\text{Opt}_{\text{SVM}}/\text{Opt}_{\text{Coupled}} - 1$
\end{itemize}
\begin{table}[h]
    \centering
    \addtolength{\tabcolsep}{1.6pt}
    \begin{tabular}{ccccccccc}
        \toprule
        \toprule
        \multicolumn{3}{c}{\textbf{Thresholds}} & & \multicolumn{5}{c}{\textbf{Metrics}}\\
        \midrule[0.05pt]
        $d_{bp, 1}$ & $d_{bp, 2}$ & $\theta_{tr}$ & & \# LS & \# QP & \# I & Wall/I (ms) & Opt. (\%) \\
        \midrule[0.05pt]
        0 & 0.15 & 5\textdegree & & \textbf{21} & 31 & 13 & 1.02 ms & 0.04\%\\
        0.15 & 0.15 & 5\textdegree & & 58 & 0.6 & \textbf{11} & \textbf{0.69 ms} & 0.34\% \\
        1 & 0.15 & 5\textdegree & & 413 & \textbf{0} & 13 & 0.75 ms & 0.31\% \\
        0 & 0 & 5\textdegree & & \textbf{21} & 1.45 & \textbf{11} & 0.82 ms & 0.33\% \\
        0 & 1 & 5\textdegree & & \textbf{21} & 271 & 14.5 & 2.65 ms & \textbf{0.00\%} \\
        0 & 1 & 0.1\textdegree & & \textbf{21} & 310 & 16 & 2.69 ms & \textbf{0.00\%} \\
        0 & 1 & 45\textdegree & & 25 & 186 & \textbf{11} & 2.37 ms & 0.31\% \\
        \midrule[0.05pt]
        \multicolumn{3}{c}{\textit{Coupled}} & & - & - & 17 & 2.4 ms & 0.00\% \\
        \bottomrule
        \bottomrule
    \end{tabular}
    \caption{Overview of the effect of different thresholds of the filter parameters on the performance aspects of solving an optimal control problem using the decoupled approach.}
    \label{tab:effects_filter}
\end{table}

The numbers in Table~\ref{tab:effects_filter} lead to various conclusions. Firstly, the effect of the broad-phase filter(s) is clearly visible. If we increase its threshold, the number of times that the LS or QP is solved increases. This could be expected since new hyperplanes for more robot-obstacle pairs must be computed, which can be derived from Algorithm~\ref{alg:compute_svm}. The subsequent effect on the wall time per iteration is two-fold, depending on whether the LS or QP is solved. The increase in \# LS has a relative small effect on the wall time per iteration, as solving a linear system of equations of six-by-six (one vertex of the robot + four corner vertices of an obstacle) only takes a few microseconds. Solving the QP however, takes tens of microseconds and can have a detrimental effect on the wall time per iteration if it needs to be solved many times. \\

Secondly, the decoupling adds little conservativeness compared to the gain in computational efficiency. The increase in optimality can be attributed to two effects: (1) the trust-region filter that prevents some hyperplanes to change close to the optimal solution and thus does not allow the solver to find a more optimal solution, and (2) the decoupled hyperplane constraints are a $0^{\text{th}}$-order approximation of the coupled constraint, removing derivative information on how a change in hyperplane(s) affects a trajectory. The second reason increases optimality since the solver cannot exploit this information in the final stages of convergence to find a better local optimum. Lastly, Table~\ref{tab:effects_filter} shows that the amount of iterations required to find a trajectory is decreased by using the decoupled approach. This can be attributed to the loss of non-convex constraints and fewer optimisation variables in the decoupled approach. For the following experiments, the filter parameters $[0.15, 0.15, 5\degree]$ are chosen. They form a good balance between the beneficial effect on wall time and number of iterations but allow for a change in hyperplanes compared to the other, more extreme thresholds. \\

\subsubsection*{Computational Effort}
\begin{figure}
    \centering
    \includegraphics[width=1.0\linewidth]{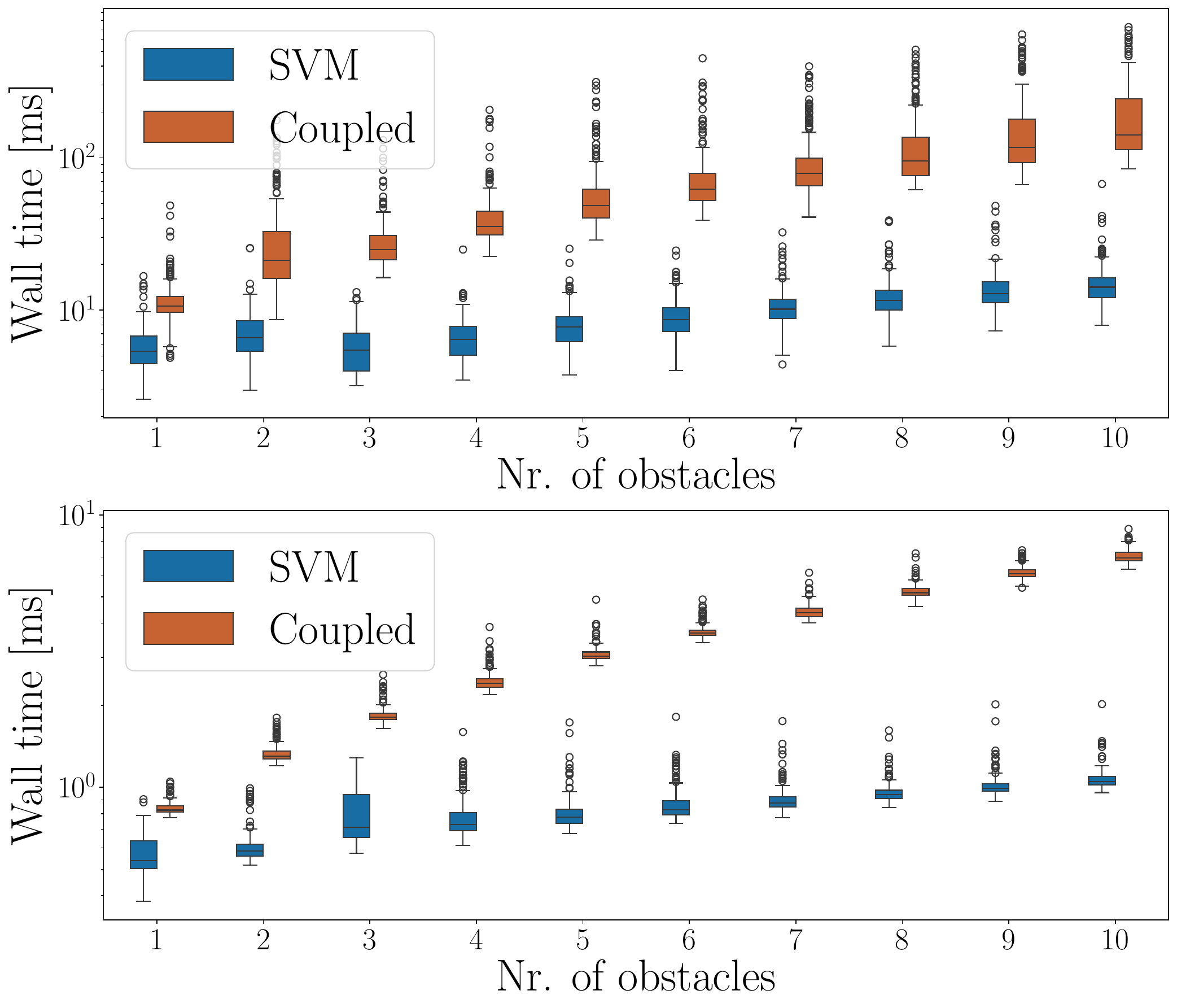}
    \caption{Total wall time (top) and wall time per iteration (bottom) required to solve the holonomic robot OCP with an increasing number of obstacles for both the coupled and decoupled approach.}
    \label{fig:2d_wall_iter_time}
\end{figure}

The top figure in Fig.~\ref{fig:2d_wall_iter_time} displays the variation in wall times necessary to find a solution per number of obstacles. Two conclusions can be derived. Firstly, by adopting the decoupled approach, the resulting wall times are reduced increasingly with the number of obstacles as compared to the coupled approach, from 50 \% up to 90 \%. This can be attributed both to a lower amount of iterations and the lower cost per iteration, as seen on the bottom figure. The decreased iteration cost is attributed to the reduced number of constraints and variables in the decoupled approach. Secondly, the variation in wall times is lower for the decoupled approach in contrast with the coupled approach. Since the bilinear hyperplane constraint is now transformed into a linear constraint, the overall non-convexity of the NLP is decreased. This makes the NLP easier to solve with less variation in the steps taken due to this decrease and thus more consistent computation times. It is concluded that
the reduction in computation time by removing the hyperplanes from the optimisation variables clearly outweighs the additional effort of finding hyperplanes through Eq.~\ref{eq:qpsvm} and~\ref{eq:lssvm}.


\subsection{Configuration-to-Point Manipulator Planning}
\label{sec:svm_robot}

To show the beneficial effects of decoupling collision avoidance on a more realistic example, the problem of configuration-to-point planning in a cluttered three-dimensional environment with a six-DOF robotic manipulator is considered. Geometrically, the robot is modelled as a set of capsules, where the hemispherical endpoints of the capsules are the data points for classification, similar to Fig.~\ref{fig:hyper_theorem}. Only the final capsule of the robot is considered here for collision avoidance for simplicity. However, the method can directly include more robot primitives, as long as their vertices are known.\\

For a single environment with five cuboids, twenty collision-free start configurations and end points are constructed over the robot's workspace, leading to a total of twenty trajectories that must be computed. The end point is converted into a final configuration using inverse kinematics, leading to a straight-line joint
space initialisation between start and end configurations. Nine out of the twenty initial trajectories are partially infeasible and allows us to test the method's capability to handle it. Similarly to the previous section, the hyperplanes are initialised by solving LS-SVM or QP-SVM along this initial trajectory.\\

\begin{figure}[h]
    \centering
    \includegraphics[width=0.75\linewidth]{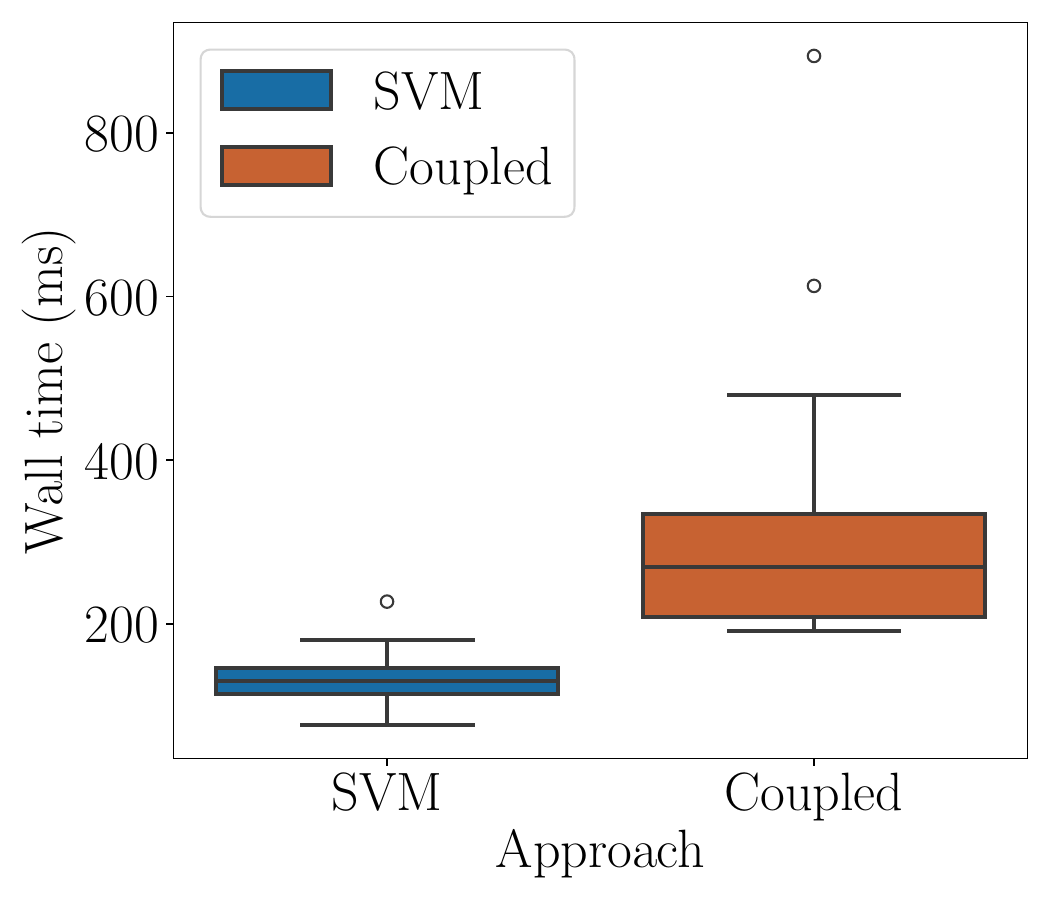}
\caption{Box plot of the wall times comparing the decoupled (SVM) and coupled approach for the robot manipulator case for twenty trajectories.}
\label{fig:robot_wall_time}
\end{figure}

Firstly, the coupled and decoupled approaches are again compared in terms of wall time required to compute a single trajectory. Fig.~\ref{fig:robot_wall_time} and Table~\ref{tab:robot_wall_time} show that the decoupled approach decreases the median computational effort by 51 \%. Table~\ref{tab:overhead_svm} shows that only 0.0 \% to 1.25 \% of this total wall time is spent on solving Eq.~\ref{eq:lssvm} with an average time of $8.71\mu s$ to solve a single eleven-by-eleven linear system. Maximally 2.55\% of the total computational time is spent on solving Eq.~\ref{eq:qpsvm}, which takes an average of $117 \mu s$ to solve. Compared to the computational efficiency of geometric toolboxes such as Coal~\cite{hppfcl} that perform similar collision checks in less than one microsecond, these numbers are large. However, considering that the total cost to find a single trajectory is in the order of several hundreds of milliseconds, it shows that the computational bottleneck is not formed by finding a solution to LS-SVM or QP-SVM and their use in these contexts is warranted. Furthermore, this approach leads to differentiable collision avoidance constraints whereas the functions from the geometric toolboxes are not differentiable and thus not applicable to gradient-based numerical NLP solvers such as Ipopt.\\

Secondly, some conservativeness is introduced by removing the hyperplanes as optimisation variables. Table~\ref{tab:robot_wall_time} shows the relative optimality of both approaches, calculated in terms of execution time of each trajectory as $T_{\text{SVM}}/T_{\text{Coupled}} - 1$ to study this possible drawback. It shows that for environments with sufficient free space, as the ones studied here, the increase in optimality is negligible and poses no drawback for using the decoupled approach.

\begin{table}[h]
    \centering
    \addtolength{\tabcolsep}{5pt}
    \begin{tabular}{rcccc}
    \toprule
    \toprule
        \textbf{Method} & Wall time & Mean & Max & Min\\
    \midrule[0.05pt]
         LS & 0.0 - 1.25 \% & 8.71 $\mu$s & 14.31 $\mu$s & 7.87 $\mu$s\\
         QP & 0.0 - 2.55 \% & 116.96 $\mu$s & 137.57 $\mu$s & 103.35 $\mu$s\\
    \bottomrule
    \bottomrule
    \end{tabular}
    \caption{Overview of the min, max and mean wall time necessary to solve Eq.~\ref{eq:qpsvm} and~\ref{eq:lssvm} and the percentage of total wall time it takes for a single trajectory.}
    \label{tab:overhead_svm}
\end{table}
\vspace{-3ex}
\begin{table}[h]
    \centering
    \addtolength{\tabcolsep}{2.5pt}
    \begin{tabular}{rcccccc}
    \toprule
    \toprule
        \textbf{Metric} & Median & LS-SVM & QP-SVM & Optimality & Success\\
    \midrule[0.05pt]
         SVM & \textbf{130.96 ms} & 0.59 ms & 0.57 ms & 0.08\% & 100\%\\
         Coupled & 269.69 ms & -  & - & \textbf{0.00\%} & 100\%\\
    \bottomrule
    \bottomrule
    \end{tabular}
    \caption{Summary of the median wall time to find a trajectory, the mean total wall time to solve Eq.~\ref{eq:qpsvm} and~\ref{eq:lssvm}, average relative optimality of the final solutions and success rates for twenty different start configurations and end poses.}
    \label{tab:robot_wall_time}
\end{table}



\vspace{-1.5ex}
\section{Conclusion and Future Work}
\label{sec:conclusion}
This paper presented an approach to reduce the computational effort involved in solving optimal control problems that use the separating hyperplane theorem for collision avoidance between convexly modelled robots and obstacles. By reformulating the theorem as a supervised classification problem, this approach allows us to remove the hyperplanes from the optimisation variables and consider them as optimisation parameters. Their value is computed using a system of linear equations and updated between iterations of a numerical NLP solver. Since only geometric information in Cartesian space is required, the approach (1) is independent of the specific robot and thus not require retraining for a new robot, (2) can easily deal with changing environments or be directly deployed in new environments without classifier update or re-training algorithms, and (3) can deal with an infeasible initialisation of OCPs.\\

The benefits of decoupling the collision avoidance from the OCP are successfully shown on two experiments. For a two-dimensional motion planning example, it successfully reduces computation times by at least 50\% for an environment with a single obstacle. For two-dimensional environments up to ten obstacles, the reduction in computation time increases up to 90\%. On a more practically relevant three-dimensional case involving a robotic manipulator, the decoupled approach reduces computation times with 51\% while introducing only a negligible decrease in optimality. On both cases, the variation in computation time is decreased. This is attributed to the substitution of non-convex collision avoidance constraints with linearised versions in the decoupled approach. \\

Future extensions to this work exist in several aspects. Firstly, the current implementation with ´lapacklu´ currently does not fully exploit the structure of Eq.~\ref{eq:lssvm}. Therefore, the wall times in Table~\ref{tab:overhead_svm} could be reduced by using a state-of-the-art linear algebra solver such as BLASFEO~\cite{Frison2018}. Secondly, the dense system of equations becomes increasingly more computationally complex to solve with a rising number of vertices. To keep the decoupling tractable, finding the hyperplanes based on the GJK algorithm is an interesting avenue to explore. Lastly, replacing the zero-order Taylor approximation of the non-convex constraint with a first- or higher-order Taylor expansion might reduce the dependency of the algorithm on the filters. A higher order approximation introduces additional gradient information of the solution of the LS-SVM system changes with respect to the trajectory, possibly improving convergence of the numerical solver.

\printbibliography






\end{document}